\documentclass[12pt,a4paper]{article}
\usepackage{amsfonts,amsmath,amssymb,indentfirst}

\usepackage{color}

\newcommand{\er}{\mathbb{R}}
\newcommand{\R}{\mathbb{R}}

\usepackage[colorlinks]{hyperref}
\hypersetup{linkcolor=blue,citecolor=blue,filecolor=black,urlcolor=blue}

\usepackage[showonlyrefs]{mathtools}
\mathtoolsset{showonlyrefs=true}

\newcommand{\mf}[1]{\mathbf{#1}}
\newcommand{\ub}{\mathbf{u}}

\newtheorem{Teorema}{Theorem}[section]
\newtheorem{Propriedade}[Teorema]{Proposition}

\usepackage[affil-it]{authblk}

\begin{document}
\title{Ground states for a nonlinear Schr\"odinger system with sublinear coupling terms}
\author{Filipe Oliveira\footnote{e-mail address: fso@fct.unl.pt}}
\affil{Mathematics Department, FCT-UNL\\Universidade Nova de Lisboa\\
Caparica Campus, 2829-516, Portugal}

\author{Hugo Tavares\footnote{email address: htavares@math.ist.utl.pt}}
\affil{CAMGSD, Mathematics Department\\
Instituto Superior T\'ecnico, Universidade de Lisboa\\
Av. Rovisco Pais, 1049-001 Lisboa, Portugal}

\maketitle

\begin{abstract}
\noindent
 We study the existence of ground states for the coupled Schr\"odinger system
 \begin{equation*}
\begin{cases}
\displaystyle -\Delta u_i+\lambda_i u_i= \mu_i |u_i|^{2q-2}u_i+\sum_{j\neq i}b_{ij} |u_j|^q|u_i|^{q-2}u_i \\
u_i\in H^1(\R^n), \quad i=1,\ldots, d,
\end{cases}
\end{equation*}
\noindent
$n\geq 1$, for $\lambda_i,\mu_i >0$, $b_{ij}=b_{ji}>0$ (the so-called ``symmetric attractive case'') and $1<q<n/(n-2)^+$. We prove the existence of a nonnegative ground state $(u_1^*,\ldots,u_d^*)$ with $u_i^*$ radially decreasing.  Moreover we show that, for $1<q<2$, such ground states are positive in all dimensions and for all values of the parameters. 
\end{abstract}

\section{Introduction}
In this paper we consider the system of $d$ equations
\begin{equation}
\label{initial}
\begin{cases}
\displaystyle -\Delta u_i+\lambda_i u_i= \mu_i |u_i|^{2q-2}u_i+\sum_{j\neq i}b_{ij}|u_j|^q|u_i|^{q-2}u_i \\
u_i\in H^1(\R^n), \quad i=1,\ldots, d
\end{cases}
\end{equation}
 with $\lambda_i, \mu_i>0$, and $b_{ij}=b_{ji}>0$, which appears in several physical contexts, namely in nonlinear optics (see for instance \cite{Ph1} and the references therein).  We also assume that
 \begin{equation}
 \label{cond}
 1<q<\left\{\begin{array}{llll}
            +\infty&if& n=1,2\\
            \\
            \displaystyle\frac{n}{n-2}&if& n\geq 3.
           \end{array}\right.
\end{equation}
Observe that the system is variational, more precisely of gradient type, as solutions can be obtain as critical points of the $C^1$--functional $$I_d:E=(H^1(\R^n))^d\to \R$$
defined by
$$
I_d(\ub)=I_d(u_1,\ldots, u_d):=\frac 12 \sum_{i=1}^d \|u_i\|_{\lambda_i}^2-\frac 1{2q} \sum_{i=1}^d  \mu_i |u_i|^{2q}_{2q}-\frac{1}{q}\mathop{\sum_{i=1}^d}_{i<j}  b_{ij}|u_iu_j|_q^q,
$$
 where $|\cdot|_{q}$ denotes the standard $L^{q}(\R^n)$ norm, while 
$$
\|u_i\|^2_{\lambda_i} :=\int_{\R^N}(|\nabla u_i|^2+\lambda_i u_i^2)\, dx \qquad i=1,\ldots, d.
$$
We will focus in the existence of \emph{ground state solutions} of \eqref{initial}, that is, solutions of the system that achieve the \emph{ground state level}
\begin{equation}\label{eq:goundstate}
c:=\inf\{ I_d(\ub):\ \ub\neq \mf{0},\ I_d'(\ub)=0\}.
\end{equation}
Another interesting question is also to know if, when $c$ is achieved, the ground state $\ub$ is \emph{nontrivial}, meaning that all its components $u_i$ are nonzero.

\medbreak

This problem has attracted a lot of attention in the last decade, specially in the particular case of $d=2$ equations:
\begin{equation}
\label{elliptic}
\left\{
\begin{array}{llll}
-\Delta u_1+\lambda_1u&=& \mu_1|u_1|^{2q-2}u_1+b|u_2|^q|u_1|^{q-2}u_1\\
-\Delta u_2+\lambda_2v&=&\mu_2 |u_2|^{2q-2}u_2+b|u_1|^q|u_1|^{q-2}u_1.
\end{array}\right.
\end{equation}
For $\mu_1=\mu_2=1$, Maia, Montefusco and Pellacci proved in \cite{MMP}  that $c$ is always achieved, while there exists a positive ground state (i.e., $u_1,u_2>0$ in $\R^n$) if $b\geq \Lambda$, for a certain $\Lambda>0$ depending on $\lambda_2/\lambda_1$. The same type of result was proved for $q=2$, $n=2,3$ by Ambrosetti and Colorado \cite{AmbrosettiColorado,AmbrosettiColoradobis}, and by de Figueiredo and Lopes for $n=1$ (see \cite{FigueiredoLopes}). On the other hand, for $q\geq 2$, there are regions where all solutions must have a null component, as it was observed for instance by Sirakov \cite{Sirakov} and Chen, Zou \cite{ChenZou}. 

\bigbreak

The optimal bounds for the existence of nontrivial ground states were found by Mandel \cite{Mandel} for every $q$ as in \eqref{cond}. More precisely, in Theorem 1 of the cited paper it is shown that there exists $\bar b:=b(\lambda_2/\lambda_1,q,n)$ such that for $b<\bar b$ all ground states have a trivial component (we will call them \emph{semitrivial}), while for $b>\bar b$ all ground states are nontrivial. For $\mu_1=\mu_2=1$ and $\lambda_2/\lambda_1=\omega^2\geq 1$, the threshold is given by the expression (see eq. (5) in \cite{Mandel})
\[
\bar b=\inf \left\{ \frac{\hat c_0^{-2q}(\|u\|^2+\|v\|_{w^2}^2)^q-|u|_{2q}^{2q}-|v|_{2q}^{2q}}{2|uv|_q^q}:\ u,v\in H^1(\R^n) \right\},
\]
where $\hat c_0:= \|u_0\| |u_0|^{-1}_{2q}$ and $u_0$ is the unique positive radially decreasing function of $-\Delta u_0+u_0=|u_0|^{2q-1}u_0$ in $\R^n$ (\cite{Unic}). It is also shown that $\bar b=0$ if $1<q<2$ (see \cite[Lemmas 1.(i) and 2-(i)]{Mandel}).

Our aim is to generalize this last result for an arbitrary number of equations. In order to state our results, let us introduce some notations.

\bigskip

We will study the minimization problem
\begin{equation}
 \label{minimization}
\inf \{I_d(u)\,:\,u\in \mathcal{N}_d\},
\end{equation}
where the so-called Nehari manifold $\mathcal{N}_d$ is defined by
$$\mathcal{N}_d:=\{\mathbf{u}\in E\,:\,\mathbf{u}\neq 0, \nabla I_d(\mathbf{u})\perp \mathbf{u}\},$$
that is, $\mathbf{u}\in \mathcal{N}_d$ if and only if $\mathbf{u}\neq 0$ and 
$$\tau_d(\mathbf{u}):=\langle  \nabla I_d(\mathbf{u}), \mathbf{u}\rangle_{L^2}=\sum_{i=1}^d \|u_i\|_{\lambda_i}^2-\Big(\sum_{i=1}^d \mu_i |u_i|_{2q}^{2q}+2\sum_{j< i}b_{ij} |u_iu_j |_q^q\Big).$$
Under condition \eqref{cond} it is classical to check that $\mathcal{N}_d$ is a manifold, and that minimizers on the Nehari manifold are ground state solutions.

\bigskip
When dealing with the system (\ref{elliptic}) it is often necessary to treat the case $n=1$ separately due to the lack of compactness of the injection $H^1_{r}(\er)\hookrightarrow L^q(\er)$, $q>2$, where $H_r^1(\er)$ denotes the space of the radially symmetric functions of  $H^1(\er)$. This lack of compactness is, in a sense, a consequence of the inequality
\begin{equation}
 \label{ineq}
|u(x)|\leq C|x|^{\frac{1-n}2}\|u\|_{H^1(\er^n)}
\end{equation}
for $u\in H_r^1(\er^n)$. Indeed, \eqref{ineq} gives no decay in the case $n=1$. However, if $u$ is also radially decreasing, it is easy to establish that
$$|u(x)|\leq C|x|^{-\frac{n}2}\|u\|_{L^2(\er^n)},$$
which provides decay in all space dimensions, hence compacity by applying the classical Strauss' compactness lemma (\cite{Strauss}). Hence, putting 
$$
H_{rd}^1(\er^n)=\{u \in H^1_r(\er^n)\,:\,u\textrm{ is radially decreasing}\},
$$ 
we get the compactness of the injection $H_{rd}^1(\er^n)\hookrightarrow L^q(\er^n)$ for \emph{all} $n\geq 1$ (see the Appendix of \cite{BL} for more details), a fact that seems not very well known. We will use this fact to present a unified approach for the problem of the energy minimization of (\ref{elliptic}), valid in all space dimensions. In fact, by putting $E_{rd}=(H_{rd}^1(\R^n))^d$ the cone of symmetric radially decreasing non-negative  functions of $E$, we will prove the following result (see also \cite{LiuWang,MMP}):


\begin{Propriedade}
\label{Teorema1}
Let $n\geq 1$ and take $q$ satisfying \eqref{cond}. There exists a minimizing sequence $(u_{1,k},\ldots,u_{d,k})\in E_{rd}$ for the minimization problem (\ref{minimization}). Furthermore, 
$(u_{1,k},\ldots,u_{d,k})\to (u_1^*,\ldots,u_d^*)\in E_{rd}$ strongly in $(H^1(\er^n))^d$. In particular
$$
I(u_1^*,\ldots,u_d^*)=\min_{\mathcal N} I(u_1,\ldots, u_d)=\min_{\mathcal{N}_d\cap E_{rd}}I(u_1,\ldots, ,u_d)=c.
$$
\end{Propriedade}
Concerning the existence of ground states with non-trivial components, we will show the following, which corresponds to our main result:

\begin{Teorema}
\label{Teorema2} 
 Let $n\geq 1$, $\lambda_i>0$, $\mu_i>0$ and $b_{ij}=b_{ji}>0$.\\ For $1<q<2$ the system \eqref{initial} admits a ground state solution $\mathbf{u}=(u_1,\cdots,u_d)\in E_{rd}$ with $u_i>0$ for all $i$. Moreover, all possible ground state solutions have nontrivial components.
\end{Teorema}

We recall that such theorem was shown by Mandel for systems with $d=2$ equations, as a consequence of a more general result, namely the characterization of the optimal threshold $\bar b$ defined before. Simplifying its approach, we will be able to prove Theorem \ref{Teorema2}, arguing by induction in the number of equations. Roughly, assuming that a subsystem of \eqref{initial} with $d-1$ equations has a certain ground state solution with nontrivial components, we will construct an element $(U_1,\ldots, U_d)\in \mathcal{N}_d$ with lower energy $I_d$.

For $d>2$ equations, the only result we are aware of for ground states (as defined in \eqref{eq:goundstate}) is the paper by Liu and Wang \cite{LiuWang}, where a ground state is proved to exist in the case $n=2,3$, $q=2$, $\lambda_1=\ldots=\lambda_m$, $\mu_1=\ldots=\mu_m$ and $\beta_{ij}=\beta$ is sufficiently large (see also \cite{Chang} for $d=3$).

We would like to observe that, for $q\geq 2$, there  are regions were the ground states are all semitrivial. Then it becomes an interesting question to study if there exist \emph{least energy nontrivial solutions} of \eqref{initial}, that is, solutions which minimize the energy among the set of all nontrivial solutions. This has been done for $d=2$ equations for instance in \cite{AmbrosettiColorado,ChenZou,LinWei,Sirakov}, and for $d\geq 2$ in \cite{SatoWang,Soave,SoaveTavares}, among others. For some recent results in this directions concerning a Schr\"odinger-KdV system, see also \cite{Colorado1,Colorado2} .


\section{Proof of Proposition \ref{Teorema1}}
We begin by observing that for $(u_1,\ldots,u_d)\in E$, $(u_1,\ldots,u_d)\neq (0,\ldots, 0)$, with $\tau_d(u_1,\ldots, u_d)\leq 0$, there exists $t\in]0,1]$ such that $(tu_1,\ldots,t u_d)\in\mathcal{N}_d$. Indeed, if $\tau_d (u_1,\ldots,u_d)=0$, we choose $t=1$. If $\tau_d(u_1,\ldots,u_d)<0$ we simply notice that 
$$
\tau(t u_1,\ldots,t u_d)=t^2\Big(\sum_{i=1}^d\|u_i\|_{\lambda_i}^2-t^{2q-2}( \sum_{i=1}^d |u_i|_{2q}^{2q}+2 \sum_{i<j} b_{ij}|fg|_q^q)\Big):=t^2T_\mf{u}(t),
$$
with $T_\mf{u}(0)>0$ and $T_\mf{u}(1)<0$.\\
Also, we notice that if $(u_1,\ldots, u_d)\in\mathcal{N}_d$,
\begin{align}\label{valori}
I(u_1,\ldots, u_d) =\Big(\frac 12-\frac 1{2q}\Big)\sum_{i=1}^d \|u_i\|^2_{\lambda_i} \\
			   =\Big(\frac 12-\frac 1{2q}\Big)\Big( \sum_{i=1}^d |u_i|_{2q}^{2q}+2 \sum_{i<j} b_{ij}|fg|_q^q)  \Big).
\end{align}
We now take a minimizing sequence $(u_{1,k},\ldots,u_{d,k})\in \mathcal{N}_d$ for the problem:
\[
\inf \{I_d(u_1,\ldots, u_d):\ (u_1,\ldots, u_d)\in \mathcal{N}_d\}.
\]
From (\ref{valori}), it is clear that this infimum is nonnegative, whence the sequence $(u_{1,k},\ldots,u_{d,k})$ is bounded in $E$.

\bigskip

We put $u_{i,k}^*$ the decreasing radial rearrangements of $|u_{i,k}|$, $i=1,\ldots, d$. It is well-known that this rearrangement preserves the $L^p$ norm ($1\leq p \leq +\infty$). Furthermore, the P\'olya-Szeg\"o inequality
$$|\nabla f^*|_2\leq |\nabla |f||_2$$ in addition with the inequality $|\nabla |f||_2\leq |\nabla f|_2$ (see \cite{Lions1}) shows that
$$
\sum_{i=1}^d \|u_{i,k}^*\|_{\lambda_i}^2\leq \sum_{i=1}^d \| u_i\|^2_{\lambda_i}.
$$
On the other hand, the Hardy-Littlewood inequality
$$\int |fg|\leq \int f^*g^*$$
combined with the monotonicity of the map $\lambda\to\lambda^q$ (see for instance \cite{Rear} for details) yields $\|fg\|_q\leq \|f^*g^*\|_q$ and, finally (as $b_{ij}>0$),
$$\tau_d(u_1^*,\ldots,u_d^*)\leq \tau_d(u_1,\ldots,u_n)=0.$$
Next, let $t_k\in]0,1]$ such that $(t_ku_{1,k}^*,\ldots,t_k u_{d,k}^*)\in \mathcal{N}_d.$
We obtain 
\begin{align*}
I(t_ku_{1,k}^*,\ldots,t_k u_{d,k}^*)&=t_k^2\Big(\frac 12-\frac 1{2q}\Big) \sum_{i=1}^d \|u_{i,k}^*\|_{\lambda_i}^2\\
						&\leq \Big(\frac 12-\frac 1{2q}\Big) \sum_{i=1}^d \|u_{i,k}\|_{\lambda_i}^2=I(u_{1,k},\ldots,u_{d,k})
\end{align*}
and we obtained a minimizing sequence $(t_ku_{1,k}^*,\ldots,t_k u_{d,k}^*)$ in $E_{rd}$, denoted again, in what follows, by $(u_{1,k},\ldots,u_{d,k})$.
Since this sequence is bounded in $E_{rd}$, up to a subsequence, $u_{i,k}\rightharpoonup u_i^*$ in $H^1(\er^n)$ weak. Also, since the injection $E_{rd}\to L^{2q}(\er^n)$ is compact, up to a subsequence, $u_{i,k}\to u_i^*$ in $L^{2q}(\er^n)$ strong, for all $n\geq 1$. Moreover, $(u_1^*,\ldots, u_d^*)\neq (0,\ldots, 0)$, since (from the definition of $\mathcal{N}_d$ and by Sobolev and Cauchy-Schwarz inequalities):
\[
\sum_{i=1}^d |u_{i,k}|^{2}_{2q}\leq C_1 \sum_{i=1}^d \|u_{i,k}\|_{\lambda_i}^2\leq C_2 \sum_{i=1}^d |u_{i,k}|_{2}^{2q}.
\]
Thus there exists $\delta>0$, independent from $k$, such that $\displaystyle\sum_{i=1}^d |u_{i,k}|_{2q} \geq \delta$, and by the strong convergence also $\displaystyle\sum_{i=1}^d |u_i^*|_{2q}\geq \delta>0$.\\
Also, since
$$
\tau(u_1^*,\ldots,u_d^*)\leq \liminf \tau(u_1,\ldots,u_d)=0,
$$
once again we can take $t\in]0,1]$ such that $(tu_1^*,\ldots,tu_d^*)\in \mathcal{N}_d$. Then,
\begin{multline*}
\inf_{\mathcal{N}_d} I_d \leq I(tu_1^*,\ldots tu_d^*)=t^2 \Big(\frac 12-\frac 1{2q}\Big) \sum_{i=1}^d  \|u_i^*\|_{\lambda_i}^2\\
\leq \Big(\frac 12-\frac 1{2q}\Big)\liminf \sum_{i=1}^d\|u_i\|_{\lambda_i}^2\leq \liminf I(u_1,\ldots, u_d)=\inf_{\mathcal{N}_d} I_d.
\end{multline*}
This implies that $(tu_1,\ldots, t u_d^*)$ is a minimizer. In particular, all inequalities above are in fact equalities: $t=1$, $(u_1^*,\ldots, u_d^*)\in \mathcal{N}_d$, $u_{i,k}\to u_i*$ in $H^1(\er^n)$ strong.\\
It is then clear that $(u_1^*,\ldots, u_d^*)$ is a ground state solution, which concludes the proof of Proposition \ref{Teorema1}. Observe also that by eventually taking the absolute value,  we can assume that $u_i^*\geq 0$, and by the strong maximum principle either $u_i>0$ or $u_i\equiv 0$. \hfill$\blacksquare$

\section{Ground states with non-trivial components. Proof of Theorem \ref{Teorema2}}

As stated in the introduction, the general result will be obtained by induction on the number of equations $d$. We begin by considering the case $d=2$. In this case, the result stated in Theorem \ref{Teorema2}  was recently obtained by Mandel in \cite{Mandel} in the case $\mu_1=\mu_2=1$. Here, we will cover this case by a different (and more direct) method, considering also arbitrary  $\mu_1,\mu_2>0$. Furthermore, as stated previously, our method will also extend easily to more general systems of $d>2$ equations.\\

\bigskip

Denote by $c_i$ the energy level of the (unique) positive ground state $u_i$ of
$$-\Delta u+\lambda_iu=\mu_i|u|^{2q-2}u.$$
Without loss of generality, we may assume that $c_1\leq c_2$.
Hence, in order to prove our result, it is sufficient to exhibit $(U_1,U_2)\in\mathcal{N}_2$, with $U_1,U_2\neq 0$ such that $I_2(U_1,U_2)<I_2(u_1,0)=c_1$.\\
For a fixed $w$ and for $\theta>0$ that will be chosen later, we begin by computing $t>0$ such that $(tu_1,t\theta w)\in \mathcal{N}_2$, that is
\begin{multline*}
\tau_2(tu_1,t\theta w)=t^2\mu_1\|u_1\|_{\lambda_1}^2+t^2\mu_2\|\theta w\|_{\lambda_2}^2\\
		-t^{2q}(\mu_1|u_1|_{2q}^{2q}+\mu_2\theta^{2q}|w|_{2q}^{2q}+2b_{12}\theta^q|u_1w|_q^q)=0,
\end{multline*}
from where we obtain that
$$t^{2q-2}=\frac{\|u_1\|_{\lambda_1}^2+\theta^2\|w\|_{\lambda_2}^2}{\mu_1|u_1|_{2q}^{2q}+\mu_2\theta^{2q}|w|_{2q}^{2q}+2b_{12}\theta^q|u_1w|_q^q}.$$
Since $u_1\in\mathcal{N}_1$, $\|u_1\|_{\lambda_1}^2=\mu_1|u_1|_{2q}^{2q}$, and we obtain
\begin{equation}
 \label{tpotencia}
t^{2q-2}=\frac{1+\theta^2C_1}{1+\mu_2\theta^{2q}C_2+2b_{12}\theta^qC_3},
\end{equation}
where 
$$C_1=\frac{\|w\|_{\lambda_2}^2}{\|u_1\|_{\lambda_1}^2},\,C_2=\frac{|w|_{2q}^{2q}}{\|u_1\|_{\lambda_1}^2}\textrm{ and } C_3=\frac{|u_1w|_q^q}{\|u_1\|_{\lambda_1}^2}.$$
Since $(tu_1,t\theta w)\in\mathcal{N}_2$, $$I(tu,t\theta w)=\Big(\frac{1}{2}-\frac 1{2q}\Big)\Big(\|tu_1\|_{\lambda_1}^2+\theta^2\|t\theta w\|_{\lambda_2}^2\Big)=t^2\Big(\frac{1}{2}-\frac 1{2q}\Big)(1+C_1\theta^2)\|u_1\|_{\lambda_1}^2,$$
and condition $I_2(tu_1,t\theta w)<I_2(u_1,0)$ is equivalent to
$$t^2\Big(\frac{1}{2}-\frac 1{2q}\Big)(1+C_1\theta^2)<1,$$
that is, in view of (\ref{tpotencia}),
$$\Big(\frac{1+\theta^2C_1}{1+\mu_2\theta^{2q}C_2+2b_{12}\theta^qC_3}\Big)^{\frac 1{q-1}}(1+C_1\theta^2)<1$$
and
$$\frac{(1+\theta^2C_1)^{q}}{1+\mu_2\theta^{2q}C_2+2b_{12}\theta^qC_3}<1.$$
Thus, we obtain
$$\frac{(1+\theta^2C_1)^{q}-1-\mu_2\theta^{2q}C_2}{\theta^q}<2b_{12}C_3.$$
By noticing that for $1<q<2$ $\displaystyle \lim_{\theta\to 0^+}\frac{(1+\theta^2C_1)^{q}-1}{\theta^q}=0$, we conclude that this condition holds for small $\theta$, which concludes the proof.\hfill$\blacksquare$

\bigskip

\noindent

We now consider system \eqref{initial} with $d>2$ equations. Given $I\subsetneq\{1,2,\ldots, d\}$ denote by $c_I$ the ground state level of the system

\[
-\Delta u_i+\lambda_i u_i= \mu_i |u_i|^{2q-2}u_i+\sum_{j\in I,j\neq i}b_{ij}|u_j|^q|u_i|^{q-2}u_i,\qquad i\in I.
\]
Let us now assume, by induction hypothesis, that there exists a ground state level $c_I$ with $\#I=d-1$ and for all $J$ with $\#J<d-1$, $c_I<c_J$. Without loss of generality, we assume that
\[
c:=c_{\{1,\ldots, d-1\}}=\min \{ c_I:\ \#I=d-1\},
\]
where $c$ is achieved by the ground state $(u_1,\ldots, u_{d-1})\in \mathcal{N}_{d-1}$, solution of
\[
-\Delta u_i+\lambda_i u_i= \mu_i |u_i|^{2q-2}u_i+\mathop{\sum_{j=1}^{d-1}}_{j\neq i} b_{ij}|u_j|^q|u_i|^{q-2}u_i,\qquad i=1,\ldots, d-1.
\]

Noticing that $I_d(u_1,\ldots,u_{d-1},0)=I_{d-1}(u_1,\ldots, u_{d-1})$, we will prove our assertion by exhibiting $(U_1,\ldots ,U_d)\in\mathcal{N}_d$, $U_i\neq 0$, such that $I_d(U_1,\ldots, U_d)<I_n(u_1,\ldots, u_{d-1},0)=c$, which guarantees that the energy level of $(U_1,\ldots, U_d)$ is inferior to the energy level of any solution of $(\ref{initial})$ with trivial components.\\
In this regard, for fixed $w\in H^1(\er^n)$, $w\neq 0$, and $\theta>0$, we choose $t>0$ such that
\begin{equation}
 \label{condicaoNehari}
(tu_1,\ldots, tu_{d-1},t\theta w)\in\mathcal{N}_d.
\end{equation}
This condition is equivalent to $\tau_d(tu_1,\ldots, tu_{d-1},t\theta w)=0$, that is
\begin{multline}
t^2\Bigl( \sum_{i=1}^{d-1}\|u_i\|_{\lambda_i}^2+\theta^2\|w\|_{\lambda_d}^2\Big)=\\
t^{2q}\Big(\sum_{i=1}^{d-1}\mu_i |u_i|_{2q}^{2q}+\mu_d\theta^{2q}|w|_{2q}^{2q}+2\mathop{\sum_{i,j=1}^{d-1}}_{j<i} b_{ij}|u_iu_j|_q^q+2\sum_{i=1}^{d-1} b_{id} \theta^q|u_iw|_q^q\Big),
\end{multline}
Since $(u_1,\ldots, u_{d-1})\in\mathcal{N}_{d-1}$, 
\[
\sum_{i=1}^{d-1}\|u_i\|_{\lambda_i}^2=\sum_{i=1}^{d-1}\mu_i |u_i|_{2q}^{2q}+2\mathop{\sum_{i,j=1}^{d-1}}_{j<i} b_{ij}|u_iu_j|_q^q,
\] 
which yields
\begin{equation}
 \label{valordet}
t^{2q-2}=\frac{1+\theta^2C_1}{1+\mu_d \theta^{2q}C_2+2\sum_{i=1}^{n-1}b_{in}\theta^q D_i},
\end{equation}
where we have put
$$C_1=\frac{\|w\|_{\lambda_n}^2}{\sum_{i=1}^{d-1}\|u_i\|_{\lambda_i}^2}, \qquad C_2=\frac{|w|_{2q}^{2q}}{\sum_{i=1}^{d-1} \|u_i\|_{\lambda_i}^2},
$$
and
$$D_i=\frac{|u_iw |_{q}^q}{\sum_{i=1}^{d-1} \|u_i\|_{\lambda_i}^2}.$$
Now, observe that, since $(tu_1,\ldots,tu_{d-1},t\theta w)\in\mathcal{N}_d$,

\begin{align*}
I_d(tu_1,\ldots,tu_{n-1},t\theta w)&=\Big(\frac{1}{2}-\frac 1{2q}\Big)\Big( \sum_{i=1}^{d-1} \|t u_i\|_{\lambda_i}^2   +\theta^2\|t\theta w\|_{\lambda_d}^2\Big)\\
&=t^2\Big(\frac{1}{2}-\frac 1{2q}\Big)\Big(1+C_1\theta^2\Big)\sum_{i=1}^{d-1} \|u_i\|_{\lambda_i}^2.
\end{align*}
Since
$$I_d(u_1,\ldots,u_{d-1},0)=I_{d-1}(u_1,\ldots,u_{d-1})=\Big(\frac{1}{2}-\frac 1{2q}\Big)\sum_{i=1}^{d-1} \| u_i\|_{\lambda_i}^2,$$
we obtain that the condition $I_d(U_1,\ldots, U_d)<c$ is equivalent to $t^2(1+C_1\theta^2)<1,$ and, in view of (\ref{valordet}), to
$$\frac{(1+\theta^2C_1)^q}{1+\mu_d \theta^{2q}C_2+2\sum_{i=1}^{d-1}b_{id}\theta^q D_i}<1$$
or
\begin{equation}
 \label{condicaofinal}
\frac{(1+\theta^2C_1)^q-1-\mu_d \theta^{2q}C_2}{\theta^q}<2\sum_{i=1}^{d-1}b_{id}D_i.
\end{equation}
which holds for $\theta$ small enough, and the proof is complete.\hfill$\blacksquare$

\bigskip

\noindent
{\bf Acknowledgment:} Filipe Oliveira was partially supported by Funda\c{c}\~ao para a Ci\^encia e Tecnologia, through contract UID/MAT/00297/2013.\\
Hugo Tavares was partially supported by Funda\c{c}\~ao para a Ci\^encia e Tecnologia through the program Investigador FCT and the project PEst-OE/EEI /LA0009/2013, as well as by the ERC Advanced Grant 2013 n. 339958 ``Complex Patterns for Strongly Interacting Dynamical Systems - COMPAT''.


\begin{thebibliography}{99}
\bibitem{Ph1} N. Akhmediev and A. Ankiewicz, Partially coherent solitons on a finite background, Phys. Rev. Lett 82, 26-61, 1999.
\bibitem{AmbrosettiColorado} A. Ambrosetti and E. Colorado, Bound and ground states of coupled nonlinear Schr\"odinger equations, C.R. Acad. Sci. Paris, Ser. 1 342, 453-458, 2006. 
\bibitem{AmbrosettiColoradobis} A. Ambrosetti and E. Colorado, Standing waves of some coupled nonlinear Schr\"odinger equations, J. London Math. Soc. 75 (1), 67-82, 2007.
\bibitem{BL} H. Berestycki and P. L. Lions, Nonlinear Scalar Field Equations, I: Existence of a Ground State, A.R.M.A. 82, 313-345, 1983.
\bibitem{Chang} J. Chang, Note on ground states of a nonlinear Schr\"odinger system, J. Math. Anal. Appl. 381, 957--962 2011.

\bibitem{ChenZou} Z. Chen and W. Zou, An optimal constant for the existence of least energy solutions of a coupled Schr\"odinger system, Calc. Var. Partial Differential Equations, 48(3-4):695-711, 2013.

\bibitem{Colorado1} E. Colorado, On the existence of bound and ground states for some coupled nonlinear Schr\"odinger-Korteweg-deVries equations, arXiv:1411.7283, 2014.
\bibitem{Colorado2} E. Colorado, Existence of Bound and Ground States for a System of Coupled Nonlinear Schr\"odinger-KdV Equations, arXiv:1410.7638, 2014.
\bibitem{FigueiredoLopes} D. de Figueiredo and O. Lopes, Solitary waves for some nonlinear Schr\"odinger systems. Ann. Inst. H. Poincar\'e Anal. Non Lin\'eaire 25(1), 149-161 (2008) 
\bibitem{Rear} H. Hardy, J. E. Littlewood and G. Polya, Inequalities, Cambridge University Press, 1952.
\bibitem{LinWeiErratum} T.-C. Lin and J. Wei, Erratum: ``{G}round state of {$N$} coupled nonlinear, {S}chr{\"o}dinger equations in {${\bf R}^n$}, {$n\leq3$}'' [{C}omm. {M}ath. {P}hys. {\bf 255} (2005), no. 3, 629--653; mr2135447], 277(2):573--576.
\bibitem{LinWei} T.-C. Lin and J. Wei,  Ground state of {$N$} coupled nonlinear {S}chr{\"o}dinger equations  in {$\bold R^n$}, {$n\leq 3$}, Comm. Math. Phys., 255(3):629--653, 2005.
 \bibitem{Lions1} P. L. Lions, The concentration-compactness principle in the calculus of variations, Part 1, Ann. Inst. H. Poincar\'e 1 (1984), 109-145. 
 \bibitem{LiuWang} Z. Liu and Z.-Q. Wang, Ground states and bound states of a nonlinear {S}chr{\"o}dinger system,  Adv. Nonlinear Stud., 10(1):175--193, 2010.
\bibitem{MMP} L.A. Maia, E.Montefusco and B. Pellacci, Positive solutions for a weakly coupled nonlinear Schr\"odinger System, J. Diff. Eq. 229, 743-767, 2006.
\bibitem{Mandel} R. Mandel, Minimal energy solutions for cooperative nonlinear Schr\"odinger systems, NoDEA - Nonlin. Differential Equations and Appl, DOI: 10.1007/s00030-014-0281-2
\bibitem{Rab} P.H. Rabinowitz, On a class of nonlinear Schr\"odinger equations, Z. Angew. Math. Phys. 43, 270-291, 1992. 
\bibitem{SatoWang} Y. Sato and Z.-Q. Wang, Least energy solutions for nonlinear {S}chr{\"o}dinger systems with  mixed attractive and repulsive couplings.
\bibitem{Sirakov} B. Sirakov, Least energy solitary waves for a system of nonlinear Schr\"dinger equations in $\R^n$, Comm. Math. Phys. 271, 199-221, 2007.
\bibitem{Soave} N. Soave, On existence and phase separation of solitary waves for nonlinear {S}chr{\"o}dinger systems modelling simultaneous cooperation and competition, Calc. Var. Partial Differential Equations, DOI:  10.1007/s00526-014-0764-3.
\bibitem{SoaveTavares} N. Soave and H. Tavares, New existence and symmetry results for least energy positive solutions of Schr\"odinger systems with mixed competition and cooperation terms,  arXiv:1412.4336.
\bibitem{Strauss} W. A. Strauss, Existence of solitary waves in higher dimensions, Comm. Math. Phys. 55, 149-162, 1977.
\bibitem{Unic} M.K. Kwong, Uniqueness of positive solutions of $-\Delta u+u=u^p$ in $\er^N$, Arch. Rat. Mech. Anal. 105, 243-266, 1989.

 \end{thebibliography}
 \end{document}